\documentclass[12pt]{amsart}
\usepackage{probpack}

\begin{document}
\title[RWRE with holding times]
       {Quenched large deviations for multidimensional 
         random walk in random environment\\ with holding times}
\author[R. FUKUSHIMA and N. KUBOTA]{Ryoki Fukushima and Naoki KUBOTA}
\address{Ryoki Fukushima\\
Department of Mathematics, Tokyo Institute of Technology, 
Tokyo 152-8551, JAPAN \\
Tel.: +81-(0)3-5734-2544\\
Fax: +81-(0)3-5734-2738\\}
\email{ryoki@math.titech.ac.jp} 
\address{Naoki Kubota\\
Department of Mathematics, Graduate School of Science and Technology, 
Nihon University, Tokyo 101-8308, JAPAN}
\email{kubota@grad.math.cst.nihon-u.ac.jp}
\thanks{The first author was supported by 
              JSPS Grant-in-Aid for Research Activity Start-up 22840019.}
\date{}


\begin{abstract}
We consider a random walk in random environment with random holding
times, that is, the random walk jumping to 
one of its nearest neighbors with some transition probability 
after a random holding time. 
Both the transition probabilities and the laws of the holding
times are randomly distributed over the integer lattice. 
Our main result is a quenched large deviation principle for the 
position of the random walk. 
The rate function is given by the Legendre transform of the
so-called Lyapunov exponents for the Laplace transform of the 
first passage time. 
By using this representation, we derive some asymptotics of the 
rate function in some special cases. 
\end{abstract}

\maketitle

\section{Introduction}

In this paper, we study large deviations for 
random walk in random environment with random holding times. 
The same problem has been studied by Dembo, Gantert, and
Zeitouni~\cite{DGZ04} in one-dimensional case. 
They assumed that the transition probabilities are uniformly elliptic
and holding times bounded away from zero but otherwise 
only quite general ergodicity and integrability conditions. 
We consider the multidimensional case with rather restrictive 
independence assumptions: the transition probability and 
holding times are i.i.d.~and mutually independent. 
On the other hand, we need a weaker ellipticity assumption
and also do not assume that the holding times are bounded below.

We now describe the setting in more detail. 
Denote by $\mathcal{P}_1$ the space of probability measures on the set 
$\{e \in \Z^d;|e|=1 \}$ of the canonical unit vectors of $\R^d$. 
Let $\Omega:=\mathcal{P}_1^{\Z^d}$ be the space equipped with the
canonical product $\sigma$-field $\mathcal{G}$ and an i.i.d.~probability
measure $\P$. 
Then, for an \textit{environment} $\omega=(\omega (x,\cdot))_{x \in \Z^d} \in \Omega$, the \textit{random walk in random environment} (\textit{RWRE} for short) is the Markov chain $(X=(X_n)_{n=0}^\infty,(P_\omega^x)_{x \in \Z^d})$ on $\Z^d$ defined as follows: 
$P_\omega^x (X_0=x)=1$ and 
\begin{align*}
 P_\omega^x (X_{n+1}=y+e|X_n=y)=\omega (y,e)
\end{align*}
for all $y \in \Z^d,\,n \in \N_0$ and $e \in \Z^d$ with $|e|=1$. 

Let $\mathcal{P}_2$ be the space of Borel probability measures on $(0,\infty )$. 
We consider the space $\Sigma:=\mathcal{P}_2^{\Z^d}$ endowed with the canonical $\sigma$-field $\mathcal{S}$ and an i.i.d.~probability measure $\bP$. 
Denote an element of $\Sigma$ by $\sigma=(\sigma_x)_{x \in \Z^d}$, and let $\tau=(\tau_n(x))_{n \in \N_0,x \in \Z^d} \in (0,\infty )^{\N_0 \times \Z^d}$ be independent random variables with $\sigma_x$ being the law of $\tau_n(x)$ for each $n \in \N_0$. 
We call $(\tau_n(x))_{n \in \N_0,x \in \Z^d}$ \textit{holding times} and 
denote by $P_\sigma^\HT$ their law, that is, 
$\sigma_x(ds)=P_\sigma^\HT(\tau_1(x) \in ds)$. 

For a random walk path $X$ and holding times 
$\tau$, we define the corresponding 
continuous-time random walk path $(Z_t)_{t \geq 0}$ as follows: 
\begin{align*}
 Z_t:=
 \begin{cases}
  X_n, &\sum_{m=0}^{n-1} \tau_m (X_m) \leq t <\sum_{m=0}^n \tau_m (X_m),\\
  \Delta, &t \geq \sum_{m=0}^\infty \tau_m (X_m),
 \end{cases}
\end{align*}
where $\sum_{m=0}^{n-1} \tau_m (X_m):=0$ if $n=0$ and $\Delta$ is the
graveyard state for $(Z_t)_{t \geq 0}$. 
This process $(Z_t)_{t \geq 0}$ is called \textit{random walk in random environment with holding times} (\textit{RWREHT} for short). 
Let $\tilde{P}_{\omega,\sigma}^x:=P_\omega^x \otimes P_\sigma^\HT$ and
denote the expectations with respect to
$\P,\,P_\omega^x,\,\bP,\,P_\sigma^\HT$ and $\tilde{P}_{\omega,\sigma}^x$
by $\E,\,E_\omega^x,\,\bE,\,E_\sigma^\HT$ and
$\tilde{E}_{\omega,\sigma}^x$, respectively. 
Throughout this paper, we make the following assumptions.
\begin{assum}\label{assumption}
\begin{enumerate}
\item $\log \min_{|e|=1} \omega (0,e) \in L^d(\P)\textrm{ and }
 \int_0^\infty s\,\sigma_0(ds) \in L^d(\bP)$,
\item The origin is in the convex hull of 
$\supp \Bigl( \law \bigl( \sum_{|e|=1} \omega (0,e)e \bigr) \Bigr)$.
\end{enumerate}
\end{assum}
By the first assumption and Jensen's inequality, we have 
\begin{align*}
 \theta_{\lambda ,\sigma}(z):
 =-\log \int_0^\infty e^{-\lambda s} \sigma_z (ds)\in L^d(\bP)
\end{align*}
for each $\lambda>0$ and $z \in \Z^d$. Note also that 
we always have $\theta_{\lambda ,\sigma}(z)>0$.
The second assumption is called {\it nestling property} and will 
be used only in the proof of the large deviation lower bound. 

We prove a large deviation principle (LDP) for the law of the scaled position
$Z_t/t$ of \textit{RWREHT} following the same strategy as in~\cite{Zer98a}. 
We introduce $H^Z(y):= \inf \{t \geq 0;Z_t=y \}$ as the first 
passage time through $y$ for the path $(Z_t)_{t \geq 0}$ and study the
asymptotics of the cumulant generating function as $y\to\infty$ first. 
Define for any $\lambda \geq 0,\,\omega \in \Omega,\, \sigma \in \Sigma$ and $x,y \in \Z^d$, 
\begin{align*}
 &e_\lambda (x,y,\omega ,\sigma )
  := \tilde{E}_{\omega ,\sigma}^x \bigl[ \exp \{ -\lambda H^Z(y) \} \1{\{ H^Z(y)<\infty \}} \bigr] ,\\
 &a_\lambda (x,y,\omega ,\sigma ):=-\log e_\lambda (x,y,\omega ,\sigma ).\\
\end{align*}

\begin{theorem} \label{thm:intro_1}
For each $\lambda \geq 0$,
there exists a nonrandom function $\alpha_\lambda:\R^d \to [0,\infty)$ such that for all $x \in \Z^d$, 
\begin{align}
\begin{split}
 \lim_{n \to \infty} \frac{1}{n} \,a_\lambda (0,nx,\omega ,\sigma )
 &= \lim_{n \to \infty} \frac{1}{n} \,\E \otimes \bE [a_\lambda (0,nx,\omega ,\sigma )]\\
 &= \inf_{n \geq 1} \frac{1}{n} \,\E \otimes \bE [a_\lambda (0,nx,\omega ,\sigma )]
  = \alpha_\lambda (x)
\end{split}
\label{eq:intro_1}
\end{align}
holds $\P \otimes \bP \hyphen \as$ and in $L^1(\P \otimes \bP)$. 
Moreover $\alpha_\lambda$ has the following properties: 
for any $q>0$ and $x,y \in \R^d$, 
\begin{align*}
 &\alpha_\lambda (qx)=q \alpha_\lambda (x),\\
 &\alpha_\lambda (x+y) \leq \alpha_\lambda (x)+\alpha_\lambda (y),
\end{align*}
and 
\begin{align*}
 |x|(-\log \bE[\exp \{ -\theta_{\lambda ,\sigma}(0) \} ])
 \leq \alpha_\lambda (x)
 \leq |x| \Bigl( \max_{|e|=1} \E[-\log \omega (0,e)]+\bE[\theta_{\lambda ,\sigma}(0)] \Bigr) .
\end{align*}
Furthermore, $\alpha_\lambda (x)$ is concave increasing in $\lambda \geq 0$ and convex in $x \in \R^d$. 
In particular, it is jointly continuous in $\lambda \geq 0$ and $x \in \R^d$. 
\end{theorem}

\begin{theorem} \label{thm:intro_2}
The following holds $\P \otimes \bP \hyphen \as$ and in $L^1(\P \otimes \bP)$: 
for all $\lambda \geq 0$ and all sequences $(x_n)_{n=1}^\infty$ of $\R^d$ with $|x_n| \to \infty$, 
\begin{align}
 \lim_{n \to \infty} \frac{a_\lambda (0,[x_n],\omega ,\sigma )-\alpha_\lambda (x_n)}{|x_n|}=0,
 \label{eq:intro_4}
\end{align}
where $[x_n]$ denotes a point in $\Z^d$ that is closest to $x_n$ 
in $l^1$-distance.
\end{theorem}

The following theorem is our main result. 
\begin{theorem} \label{thm:intro_3}
The law of $Z_t/t$ obeys the following large deviation principle (LDP)
with rate function 
\begin{align*}
 I(x):=\sup_{\lambda \geq 0} (\alpha_\lambda (x)-\lambda ):
\end{align*}
\begin{itemize}
 \item Upper bound: for any closed subset $A \subset \R^d$, we have $\P \otimes \bP \hyphen \as$, 
       \begin{align}
        \limsup_{t \to \infty} \frac{1}{t} \log \tilde{P}_{\omega ,\sigma}^0 (Z_t \in tA) \leq -\inf_{x \in A} I(x).
        \label{eq:intro_5}
       \end{align}
 \item Lower bound: for any open subset $B \subset \R^d$, we have $\P \otimes \bP \hyphen \as$, 
       \begin{align}
        \liminf_{t \to \infty} \frac{1}{t} \log \tilde{P}_{\omega ,\sigma}^0 (Z_t \in tB) \geq -\inf_{x \in B} I(x).
        \label{eq:intro_6}
       \end{align}
\end{itemize}
\end{theorem}

\subsection{Comments on the proof}
We basically follow the strategy taken in~\cite{Zer98a,Zer98b}. 
The second reference~\cite{Zer98b} is the first paper studying
quenched large deviations for multidimensional RWRE, where 
nearest neighbor nestling 
walks in i.i.d.~random environment are considered. 
After that, several generalizations were discussed by different 
methods. Varadhan~\cite{Var03} generalized Zerner's result to 
general ergodic case with bounded step size but under the uniform 
ellipticity assumption. His method is based on an application of
the subadditive theorem to certain hitting probabilities. 
Rosenbluth~\cite{Ros06} weakened the uniform ellipticity assumption 
and also obtained a variational formula for the
rate function by using a stochastic homogenization approach
originally applied to a diffusion with a random drift in~\cite{KRV06}. 
Later Yilmaz~\cite{Yil09} and Rassoul-Agha and 
Sepp\"{a}l\"{a}inen~\cite{RAS11} extended this method to establish
level 2 and level 3 LDPs, respectively.
Although our method requires rather restrictive independence assumptions 
and only proves a level 1 LDP, it has an
advantage of giving a relatively simple representation of the
rate function. 
This for instance allows us to determine the asymptotics of 
the rate function as $x\to\infty$ and $x\to 0$ in some 
special cases, see Section~\ref{rate function}.

Next, we explain the outline of our proof. 
To prove an LDP for random walk in random environment, 
it is standard to consider the Laplace transform of $H^Z(y)$. 
Indeed, the large deviation upper bound is almost immediate 
from Theorem~\ref{thm:intro_2} since for a compact set
$K\subset\R^d$,
\begin{equation*}
\begin{split}
 \tilde{P}_{\omega ,\sigma}^0(Z_t\in tK)
 &\le \#(tK\cap\Z^d) \max_{y\in tK\cap\Z^d}
 \tilde{P}_{\omega ,\sigma}^0(H^Z(y)\le t)\\
 &\le \#(tK\cap\Z^d) \max_{y\in tK\cap\Z^d}
 e^{\lambda t}e_{\lambda}(0,y,\omega,\sigma)\\
 &=\max_{y\in tK\cap\Z^d}\exp\{-t(\alpha_{\lambda}(y/t)-\lambda)+o(t)\}.
\end{split}
\end{equation*}
In order to extend this to general closed sets, we have to check 
what is called exponential tightness but it is not hard (see page~\pageref{eq:ldpupper_1}). 
The proof of Theorem~\ref{thm:intro_2} itself is based on the
fact that our $e_{\lambda}$ is the survival probability for a 
crossing RWRE with random potential, 
see~\eqref{eq:lyapunov_6}. Given this interpretation, 
one can prove Theorem~\ref{thm:intro_2} similarly to 
those in~\cite{Zer98a,Zer98b}.

The proof of the lower bound is a bit more complicated. 
The key to the proof is that, whenever $\alpha_{\lambda}(y)$ is 
differentiable in $\lambda$, we have
\begin{equation*}
\tilde{P}_{\omega ,\sigma}^0\left(H^Z(ty)/t \in 
\frac{\partial}{\partial\lambda}\alpha_{\lambda}(y)({1-\epsilon}, 1)
\right)
=\exp\left\{-t \left(\alpha_{\lambda}(y)-
\lambda\frac{\partial}{\partial\lambda}\alpha_{\lambda}(y)\right)+o(t)\right\} 
\end{equation*}
for any $\epsilon>0$, where $H^Z(ty)$ denotes the first time for $Z_{\cdot}$
to hit $[ty]$ (this easily follows from Lemma~\ref{lem:ldplower_0} below). 
This means that we know the cost for the random walk to make 
a crossing in the speed 
$1/\frac{\partial}{\partial\lambda}\alpha_{\lambda}(y)$.
In particular, if $\frac{\partial}{\partial\lambda}\alpha_{\lambda}(y)=1$ 
for some $\lambda$, then the above is already quite
close to the lower bound since such a $\lambda$ maximizes 
$\alpha_{\lambda}(y)-\lambda$. 
However we do not know if this is the case in general since
\begin{enumerate}
 \item[(1)] it is hard to check the differentiability of
       $\alpha_{\lambda}(y)$ and 
 \item[(2)] it may happen that
$\alpha'_{\lambda+}(y)|_{\lambda=0}=
\sup_{\lambda\ge 0}\alpha'_{\lambda+}(y)<1$. 
\end{enumerate}
To circumvent the first issue, we choose two differentiability 
points slightly above and slightly below a maximizer of 
$\lambda \mapsto \alpha_{\lambda}(y)-\lambda$. 
Then, we let the walker move slower than needed toward an 
intermediate point and faster on the rest of the way to $[ty]$ 
to achieve the expected speed, see the proof of~Lemma~\ref{lem:ldplower_2}.
As for the second issue, we find a {\em trap} around $[ty]$, 
that is, a region where the walker can spend time with 
relatively high probability. This is the content of 
Lemma~\ref{lem:ldplower_3} and
we need the nestling assumption precisely here. 

\subsection{Notation}
For $x=(x_1,\dots,x_d)$ in $\R^d$, we write $|x|:=|x_1|+\dots+|x_d|$. 
Open $l_1$-balls with center $x \in \R^d$ and radius $r \geq 0$ are
denoted by $B(x,r)$ and closed balls by $\bar{B}(x,r)$. 
We write $[x]$ for a lattice site with minimal $l_1$-distance from
$x$ chosen by some deterministic rule. Note that always $|x-[x]| \leq d/2$. 
Similarly, let $[A]:=\{ [x];x \in A \}$ for each subset $A \subset \R^d$. 
\section{Lyapunov exponents}
In this section, we show Theorem~\ref{thm:intro_1}. 
We start with the triangle inequality and integrability properties for $a_\lambda$. 
To this end, let $H^X(y):= \inf \{n \geq 0;X_n=y \}$ be the first passage time through $y$ for the random walk $(X_n)_{n=0}^\infty$. 
Then
\begin{align*}
 H^Z(y)=\sum_{m=0}^{H^X(y)-1} \tau_m (X_m)
\end{align*}
holds on the event $\{ H^X(y)<\infty \}=\{ H^Z(y)<\infty \}$ 
and hence by Fubini's theorem, 
\begin{align}
 e_\lambda (x,y,\omega ,\sigma )
 = E_\omega^x \Biggl[ \exp \Biggl\{ -\sum_{m=0}^{H^X(y)-1} \theta_{\lambda ,\sigma} (X_m) \Biggr\}
   \1{\{ H^X(y)<\infty \}} \Biggr] .
 \label{eq:lyapunov_6}
\end{align}
\begin{lemma} \label{lem:lyapunov_1}
For any $\lambda \geq 0,\,x,y,z \in \Z^d,\,\omega \in \Omega$ and $\sigma \in \Sigma$, 
\begin{align}
 a_\lambda (x,y,\omega ,\sigma ) \leq a_\lambda (x,z,\omega ,\sigma )+a_\lambda (z,y,\omega ,\sigma ).
 \label{eq:lyapunov_1}
\end{align}
Moreover, if $d=1$ and $x \leq z \leq y$ or $y \leq z \leq x$, then equality holds in \eqref{eq:lyapunov_1}. 
\end{lemma}
\begin{proof}
Set $H_z^X(y):=\inf \{ n \geq H^X(z);X_n=y \}$. 
Using the fact that $H_z^X(y)\ge H^X(z)$ and the strong Markov property, 
we have 
\begin{align}
\begin{split}
 e_\lambda (x,y,\omega ,\sigma )
 &\geq E_\omega^x \Biggl[ \exp \Biggl\{ -\sum_{n=0}^{H_z^X(y)-1} \theta_{\lambda ,\sigma}(X_n) \Biggr\}
                          \1{\{ H_z^X(y)<\infty \}} \Biggr] \\
 &= e_\lambda (x,z,\omega ,\sigma ) \,e_\lambda (z,y,\omega ,\sigma ).
\end{split}
\label{eq:lyapunov_2}
\end{align}
By taking logarithm, this proves \eqref{eq:lyapunov_1}. 
If $d=1$ and $x \leq z \leq y$ or $y \leq z \leq x$, then equality holds
 in \eqref{eq:lyapunov_2} since the random walk $(X_n)_{n=0}^\infty$ has
 to go through $z$ before reaching $y$. 
\end{proof}

\begin{lemma} \label{lem:lyapunov_2}
Let $\lambda \geq 0$ and $p \geq 1$. 
Then $a_\lambda (0,x,\omega,\sigma) \in L^d(\P \otimes \bP)$ holds for $x \in \Z^d$. 
Moreover, the collection of random variables $a_\lambda (0,x,\sigma,\omega)/|x|,\,x \in \Z^d \setminus \{ 0 \}$ is uniformly integrable under $\P \otimes \bP$ and we have for all $x \in \Z^d$, 
\begin{align}
 c_1(\lambda )|x|
 \leq \E \otimes \bE[a_\lambda (0,x,\omega ,\sigma )]
  \leq c_2(\lambda )|x|,
\label{eq:lyapunov_3}
\end{align}
where 
\begin{align*}
 &c_1(\lambda ):=-\log \bE[\exp \{ -\theta_{\lambda ,\sigma} (0) \}] ,\\
 &c_2(\lambda ):=\max_{|e|=1} \E[-\log \omega (0,e)]+\bE[\theta_{\lambda ,\sigma}(0)].
\end{align*}
\end{lemma}
\begin{proof}
Let $x \in \Z^d \setminus \{ 0 \}$. 
By forcing the walker to follow a nearest neighbor path $(0=r_0,r_1,\dots,r_m=x)$ from $0$ to $x$ with minimal length $m=|x|$, we have 
\begin{align*}
 e_\lambda (0,x,\omega ,\sigma )
 \geq \exp \biggl\{ -\sum_{n=0}^{m-1} \theta_{\lambda ,\sigma}(r_n) \biggr\} \prod_{n=0}^{m-1} \omega (r_n,r_{n+1}-r_n).
\end{align*}
It follows that 
\begin{align}
 \frac{a_\lambda (0,x,\omega ,\sigma )}{|x|}
 \leq -\frac{1}{m} \sum_{n=0}^{m-1} \log \omega (r_n,r_{n+1}-r_n)+\frac{1}{m} \sum_{n=0}^{m-1} \theta_{\lambda ,\sigma}(r_n)
 \label{eq:lyapunov_4}
\end{align}
and hence $a_\lambda (0,x,\omega ,\sigma) \in L^d(\P \otimes \bP)$ 
by Assumption~\ref{assumption}-(1).
Moreover, Jensen's inequality implies that for any $\gamma \geq 0$, 
\begin{align}
\begin{split}
 &\E \otimes \bE \Biggl[ \Biggl( -\frac{1}{m} \sum_{n=0}^{m-1} \log \omega (r_n,r_{n+1}-r_n)
                                 +\frac{1}{m} \sum_{n=0}^{m-1} \theta_{\lambda ,\sigma}(r_n)-\gamma \Biggr)_+ \Biggr] \\
 &\leq \frac{1}{m} \sum_{n=0}^{m-1} \E \otimes \bE [(-\log \omega (r_n,r_{n+1}-r_n)
       +\theta_{\lambda ,\sigma}(r_n)-\gamma )_+]\\
 &\leq \max_{|e|=1} \E \otimes \bE [(-\log \omega (0,e)+\theta_{\lambda ,\sigma}(0)-\gamma )_+].
\end{split}
\label{eq:lyapunov_5}
\end{align}
From \eqref{eq:lyapunov_4} and \eqref{eq:lyapunov_5} with $\gamma=0$, 
the right inequality of~\eqref{eq:lyapunov_3} follows. 
Let us proceed to the proof of the uniform integrability. 
Note first that for any $\gamma \geq 0$, 
\begin{align*}
 &\limsup_{M \to \infty} \sup_{x \in \Z^d \setminus \{ 0 \}}
  \E \otimes \bE \Biggl[ \frac{a_\lambda (0,x,\omega ,\sigma )}{|x|}
                         \1{\{ a_\lambda (0,x,\omega ,\sigma )/|x|>M \}} \Biggr] \\
 &\leq \limsup_{M \to \infty} \sup_{x \in \Z^d \setminus \{ 0 \}}
       \Biggl\{ \gamma \,\P \otimes \bP \Biggl( \frac{a_\lambda (0,x,\omega ,\sigma )}{|x|}>M \Biggr) \\
 &\qquad \qquad \qquad \qquad \ \ 
                +\E \otimes \bE \Biggl[ \Biggl( \frac{a_\lambda (0,x,\omega ,\sigma )}{|x|}
                                                -\gamma \Biggr)_+ \Biggr] \Biggr\}.
\end{align*}
By estimating the first term using Markov's inequality together with 
the right inequality of~\eqref{eq:lyapunov_3} and the second term 
using \eqref{eq:lyapunov_4} and \eqref{eq:lyapunov_5}, 
we see that the above right hand side is smaller than
\begin{align*}
 &\limsup_{M \to \infty} \frac{\gamma}{M} \,c_2(\lambda )
       +\max_{|e|=1} \E \otimes \bE [(-\log \omega (0,e)+\theta_{\lambda ,\sigma}(0)-\gamma )_+]\\
 &= \max_{|e|=1} \E \otimes \bE [(-\log \omega (0,e)+\theta_{\lambda ,\sigma}(0)-\gamma )_+].
\end{align*}
Since $\log \min_{|e|=1} \omega (0,e) \in L^1(\P)$ and 
$\theta_{\lambda,\sigma} (0) \in L^1(\bP)$, 
$\E \otimes \bE [(-\log \omega (0,e)+\theta_{\lambda ,\sigma}(0)-\gamma)_+]$ 
tends to zero as $\gamma \to \infty$ for any
 $e \in \Z^d$ with $|e|=1$ by Lebesgue's dominated convergence theorem. 
We thereby find that the collection of random variables $a_\lambda (0,x,\omega ,\sigma)/|x|,\,x \in \Z^d \setminus \{ 0 \}$ is uniformly integrable under $\P \otimes \bP$. 

Finally, we show the left inequality of~\eqref{eq:lyapunov_3}. 
We introduce for $\lambda \geq 0,\,\omega \in \Omega,\,\sigma \in 
\Sigma$ and $x,y \in \Z^d$ the path measure 
\begin{align*}
 \hat{P}_{\lambda ,\omega ,\sigma}^{x,y}(dX_\cdot )
 := e_\lambda (x,y,\omega ,\sigma )^{-1} \exp \Biggl\{ -\sum_{n=0}^{H^X(y)-1} \theta_{\lambda ,\sigma}(X_n) \Biggr\}
    \1{\{ H^X(y)<\infty \}} P_\omega^x (dX_\cdot )
\end{align*}
and denote its expectation by $\hat{E}_{\lambda ,\omega ,\sigma}^{x,y}$. 
In addition, let us define
\begin{align*}
 \mathcal{A}(x,X_\cdot):=\{X_0, \ldots, X_{H^X(x)}\}.
\end{align*}
Since $\# \mathcal{A}(x,X_\cdot) \geq |x|$ $P_\omega^0 \hyphen \as$, we get from Jensen's inequality that
\begin{align*}
 &c_1(\lambda )|x|\\
 &\leq \E \otimes \bE \Bigl[ \log \hat{E}_{\lambda ,\omega ,\sigma}^{0,x} [ \exp \{ c_1(\lambda )
                             \,\# \mathcal{A}(x,X_\cdot ) \} ] \Bigr] \\
 &\leq \E \otimes \bE [a_\lambda (0,x,\omega ,\sigma )]\\
&\quad +\log\E \otimes \bE \otimes E_\omega^0\Biggl[ 
\exp\Biggl\{c_1(\lambda )\,\# \mathcal{A}(x,X_\cdot )
-\sum_{n=0}^{H^X(y)-1} \theta_{\lambda ,\sigma}(X_n) \Biggr\} 
    \1{\{ H^X(y)<\infty \}}\Biggr].
\end{align*}
Clearly $\sum_{n=0}^{H^X(y)-1} \theta_{\lambda ,\sigma}(X_n)\ge
\sum_{y\in \mathcal{A}(x,X_\cdot )} \theta_{\lambda ,\sigma}(y)$.
From this and the independence of 
$\{\theta_{\lambda ,\sigma}(y)\}_{y\in\Z^d}$ under $\bP$, 
the above right-hand side is bounded by
\begin{align*}
 &\E \otimes \bE [a_\lambda (0,x,\omega ,\sigma )]
       +\log \E \otimes E_\omega^0 \Biggl[ \prod_{y \in \mathcal{A}(x,X_\cdot )}
                                    \bE [\exp \{ c_1(\lambda )-\theta_{\lambda ,\sigma}(y) \} ] \Biggr] \\
 &= \E \otimes \bE [a_\lambda (0,x,\omega ,\sigma )],
\end{align*}
where the last equality is due to the choice of $c_1(\lambda)$.
This proves the left inequality of \eqref{eq:lyapunov_3}. 
\end{proof}

Now we are in a position to prove Theorem \ref{thm:intro_1}.
\begin{proof}[\bf Proof of Theorem \ref{thm:intro_1}]
Given Lemma~\ref{lem:lyapunov_1} and \ref{lem:lyapunov_2}, the proof
 goes along the same lines as that of \cite[Proposition 4]{Zer98a} or
 \cite[Proposition 3]{Zer98b}. 
Namely, we first prove~\eqref{eq:intro_1} using the subadditive 
ergodic theorem. Then $\alpha_{\lambda}(qx)=q\alpha_{\lambda}(x)$ 
follows for $q\in\N$ and $x\in\Z^d$ by stationarity.
Finally, we extend $\alpha_{\lambda}(\cdot)$ to $\Q^d$ by 
$\alpha_{\lambda}(x/q)=\alpha_{\lambda}(x)/q$ and then to $\R^d$ by continuity. 
The convexity of $\alpha_{\lambda}(\cdot)$ follows from 
Lemma~\ref{lem:lyapunov_1} and 
the properties of $\alpha_{\lambda}$ as a function of $\lambda$ follow
from those of $a_{\lambda}$.
See the aforementioned references for details.
\end{proof}
\section{Shape theorem}
Our goal in this section is to prove Theorem \ref{thm:intro_2}, which is called the shape theorem, and to derive its generalizations. 
To this end, we recall the next lemma which plays the role of the maximal lemmas for random walk in a nonnegative random potential and random walk in random environment. 
Let us define the random distance
\begin{equation*}
 d_\lambda (x,y,\omega ,\sigma ):=
 \max \{ a_\lambda (x,y,\omega ,\sigma ),a_\lambda (y,x,\omega ,\sigma ) \}
\end{equation*}
for $x,y\in\Z^d$.
\begin{lemma} \label{lem:shape_1}
For each $\lambda \geq 0$, there is a positive constant $c_3(\lambda)$ such that the following holds $\P \otimes \bP \hyphen \as$: 
for any $\epsilon \in \Q \cap (0,1)$ there exists a positive number $R=R(\lambda,\epsilon,\omega,\sigma)$ such that 
\begin{align}
 \sup \{ d_\lambda ([x],[y],\omega ,\sigma );y \in \R^d ,\,|x-y| \leq \epsilon |x| \}
 <c_3(\lambda )\epsilon |x|
 \label{eq:shape_1}
\end{align}
holds for all $x \in \R^d$ with $|x|>R$. 
\end{lemma}
\begin{proof}
Since the proof is the same as that of \cite[Lemma 7]{Zer98a} or
\cite[Lemma 6]{Zer98b}, we only give an outline. 
Let us first consider the case $d\ge 2$. In this case for each 
$x, y \in \Z^d$, $x\neq y$, there are $2d$ self-avoiding paths
$x^{(i)}=(x^{(i)}_0=x, \ldots, x^{(i)}_{m_i}=y)$
connecting $x$ and $y$ which contain less than $|x-y|+9$ points
and mutually disjoint except for the starting and end points $x$ and $y$
(see~\cite{Kes86}, p.135).
By forcing the random walk to follow one of these paths, we get
\begin{equation*}
e_\lambda(x,y,\omega,\sigma)
\ge \max_i\prod_{k=0}^{m_i-1}\omega(x^{(i)}_k,x^{(i)}_{k+1})
e^{-\theta_{\lambda,\sigma}(x^{(i)}_k)}. 
\end{equation*}
Since the random variables appearing on the right-hand side are 
independent for different $(k,i)$'s, Chebyshev's inequality yields
\begin{equation}
\begin{split}
 &\P \otimes \bP\left(a_\lambda(x,y,\omega,\sigma)>2c(|x-y|+8)\right)\\
 &\quad\le
 \left(\frac{(|x-y|+8)\max_{|z|=1}{\rm Var}(-\log \omega(0,z)
 +\theta_{\lambda,\sigma}(0))}{c^2(|x-y|+8)^2}\right)^{2d}\\
 &\quad=O(|x-y|^{-2d})\label{lem:shape_1_proof}
\end{split}
\end{equation}
as $|x-y|\to\infty$, where
\begin{equation*}
 c=\E\left[-\log \min_{|z|=1}\omega(0,z)\right] 
 +\bE[\theta_{\lambda,\sigma}(0)].
\end{equation*}
Next, fix a finite subset $Z\subset \{z\in\R^d: |z|=1\}$ such that
$\bigcup_{x\in Z}B(x,\epsilon)\supset \{z\in\R^d: |z|=1\}$ and for
$x\in \R^d$, define $Y_x=\{y\in\Z^d:|x-y|<3\epsilon|x|\}$.
Then, since $\#\{(x,y):x\in nZ, y\in Y_x\}=O(n^d)$,
~\eqref{lem:shape_1_proof} and 
the Borel-Cantelli lemma imply that 
$\P\otimes\bP$-a.s., there exists an $n_0\in\N$ such that 
for all $n\ge n_0$, $x\in nZ$ and $y\in Y_x$, 
\begin{equation*}
 a_\lambda(x,y,\omega,\sigma)\le 2c(|x-y|+8)
 \le 6c\epsilon |x|+16c.
\end{equation*}
Finally, this is extended to all $x,y$ with $|x-y|\le \epsilon |x|$
with the help of the triangle inequality and the fact that
$\bigcup_{n\in \N}\bigcup_{x\in nZ}Y_x$ covers the whole space 
except for a bounded set. 

When $d=1$, the above argument does not work since for any 
$x,y \in \Z$ there is only one self-avoiding path. 
On the other hand since $a_\lambda$ becomes additive, if we have
\begin{equation*}
 a_\lambda(x,[(1+\epsilon)x],\omega,\sigma)\le {\rm const.}\epsilon|x|, 
\end{equation*}
then the same upper bound holds for $a_\lambda(x,y,\omega,\sigma)$ 
with $y$ lying between $x$ and $(1+\epsilon)x$.
We introduce a geometrically growing sequence defined by 
$x_0=0$, $x_{n+1}=[(1+\epsilon)x_n]+1$.
Then, by using a simple upper bound
\begin{equation*}
 a_\lambda(x_n,x_{n+1},\omega,\sigma)
\le \sum_{k=x_n}^{x_{n+1}}\left[-\log\omega(k,k+1)
+\theta_{\lambda,\sigma}(k)\right]
\end{equation*}
and a classical result for a sum of i.i.d.~random variables 
(see, e.g.~Section~6.8.5 of~\cite{Pet95}), 
one can deduce that $\P\otimes \bP$ almost surely, 
$a_\lambda(x_n,x_{n+1},\omega,\sigma)\le {\rm const.}|x_n-x_{n+1}|$
holds for all sufficiently large $n$.
This implies $a_\lambda(x_n,y,\omega,\sigma)\le {\rm const.}\epsilon|x_n|$
for all large $n$ and $y\in[x_n, x_{n+1}]$, as explained above,
and then it can be extended to any $x,y\in\Z$ with 
$|x-y|\le \epsilon|x|$ again by the triangle inequality. 
\end{proof}
\begin{proof}[\bf Proof of Theorem \ref{thm:intro_2}]
Given Lemma \ref{lem:shape_1}, one can prove the $\P \otimes \bP \hyphen \as$ convergence by the same strategy as in \cite[Theorem A]{Zer98a}. 
Then the $L^1(\P \otimes \bP)$-convergence follows from the $\P \otimes \bP \hyphen \as$ convergence by uniform integrability provided by Lemma \ref{lem:lyapunov_2}. 
\end{proof}

We next consider a generalization of Theorem \ref{thm:intro_2} for point-to-set distances instead of point-to-point distances. 
Let us define $e_\lambda (x,K,\omega,\sigma)$ for a nonempty subset
$K\subset\R^d$ as in \eqref{eq:lyapunov_6} but with $H^X(y)$ replaced 
by $H^X(K):=\inf \{ H^X(y);y \in K \}$. 
Furthermore, we write $a_\lambda (x,K,\omega,\sigma)$ for $-\log e_\lambda (x,K,\omega,\sigma)$ and denote the distance between $x$ and $K$ by $\dist (x,K):=\inf \{ |x-y|;y \in K \}$. 
Given Theorem \ref{thm:intro_2}, one can prove the following corollary by the same way as~\cite[Corollary 16]{Zer98a}. 

\begin{corollary} \label{cor:shape_1}
Let $\lambda \geq 0$ and $(K_n)_{n=1}^\infty$ be a sequence of subsets 
of $\R^d$ such that $K_n \not= \emptyset$ and $\dist (0,K_n) \to \infty$. 
Then, 
\begin{align*}
 \lim_{n \to \infty} \frac{a_\lambda (0,[K_n],\omega ,\sigma )-\inf_{x \in K_n} \alpha_\lambda (x)}{\dist (0,K_n)}
 = 0
\end{align*}
$\P \otimes \bP \hyphen \as$ 
\end{corollary}

Let us finally extend Theorem~\ref{thm:intro_2} to a directionally 
uniform version. 
The shape theorem will be used to relate crossing costs to
the Lyapunov exponent in the proof of the large deviation
lower bound. 
However, Theorem~\ref{thm:intro_2} does not suffice as it is. 
As we explained in the introduction, we shall divide the 
crossing into two pieces and for the second piece, 
we need a shape theorem with {\em moving starting points}. 

\begin{corollary} \label{cor:shape_2}
Let $x \in \Q^d \setminus \{ 0 \}$ and suppose that $\rho_1,\rho_2 \in \R$ satisfy $0 \leq \rho_1<\rho_2$. 
Then, the following holds $\P \otimes \bP \hyphen \as$: 
for all $\lambda \geq 0$ and all sequences $(x_n)_{n=1}^\infty$ of $\R^d$ with $x_n/n \to x$, 
\begin{align*}
 \lim_{n \to \infty}
 \frac{a_\lambda ([\rho_1 x_n],[\rho_2 x_n],\omega ,\sigma )-(\rho_2 -\rho_1 )\alpha_\lambda (x_n)}{|x_n|}=0.
\end{align*}
\end{corollary}
\begin{proof}
By the continuity and the homogeneity of $\alpha_{\lambda}(\cdot)$, 
it suffices to prove
that $\P \otimes \bP \hyphen$almost surely, 
\begin{align}
 \lim_{n \to \infty}
 \frac{a_\lambda ([\rho_1 x_n],[\rho_2 x_n],\omega ,\sigma )}{n}=(\rho_2 -\rho_1 )\alpha_\lambda (x)
 \label{eq:shape_5}
\end{align}
holds for all $\lambda \geq 0$ and all sequences $(x_n)_{n=1}^\infty$ 
of $\R^d$ with $x_n/n \to x$. 
Thanks to Lemma~\ref{lem:lyapunov_1} and Theorem~\ref{thm:intro_2}, we
know that the lower bound 
\begin{align*}
 &\liminf_{n \to \infty} \frac{1}{n} a_\lambda ([\rho_1 x_n],[\rho_2 x_n],\omega ,\sigma )\\
 &\geq \liminf_{n \to \infty} \frac{1}{n}
       \bigl( a_\lambda (0,[\rho_2 x_n],\omega ,\sigma )-a_\lambda (0,[\rho_1 x_n],\omega ,\sigma ) \bigr) \\
 &= (\rho_2 -\rho_1 )\alpha_\lambda (x)
\end{align*}
is valid. 
To show the upper bound
\begin{align}
 \limsup_{n \to \infty} \frac{1}{n} a_\lambda ([\rho_1 x_n],[\rho_2 x_n],\omega ,\sigma )
 \leq (\rho_2 -\rho_1 )\alpha_\lambda (x),
 \label{eq:shape_2}
\end{align}
note that we have for $K \in \N$ with $Kx \in \Z^d$, 
\begin{align}
\begin{split}
 a_\lambda ([\rho_1 nx],[\rho_2 nx],\omega ,\sigma )
 &\leq \sum_{m=\lceil \rho_1 n/K \rceil}^{\lfloor \rho_2 n/K \rfloor -1} a_\lambda (mKx,(m+1)Kx,\omega ,\sigma )\\
 &\quad
       +a_\lambda ([\rho_1 nx],\lceil \rho_1 n/K \rceil Kx,\omega ,\sigma )\\
 &\quad
       +a_\lambda (\lfloor \rho_2 n/K \rfloor Kx,[\rho_2 nx],\omega ,\sigma ).
\end{split}
\label{eq:shape_3}
\end{align}
by Lemma~\ref{lem:lyapunov_1}. 
Birkhoff's ergodic theorem shows that 
\begin{align*}
 \lim_{n \to \infty} \frac{1}{n} \sum_{m=\lceil \rho_1 n/K \rceil}^{\lfloor \rho_2 n/K \rfloor -1}
 a_\lambda (mKx,(m+1)Kx,\omega ,\sigma )
 = \frac{\rho_2 -\rho_1}{K} \,\E \otimes \bE [a_\lambda (0,Kx,\omega ,\sigma )]
\end{align*}
holds $\P \otimes \bP \hyphen \as$ 
On the other hand, we know that for any $\epsilon \in \Q \cap (0,1)$ 
and sufficiently large $n$, 
\begin{align*}
 \biggl| [\rho_1 nx]-\biggl\lceil \frac{\rho_1 n}{K} \biggr\rceil Kx \biggr|
 &\leq |[\rho_1 nx]-\rho_1 nx|+\biggl| \frac{\rho_1 n}{K}-\biggl\lceil \frac{\rho_1 n}{K} \biggr\rceil \biggr| K|x|\\
 &\leq \frac{d}{2}+K|x| \leq \epsilon |[\rho_1 nx]|
\end{align*}
and the same estimate with $\rho_1$ replaced by $\rho_2$ and 
the ceiling function by the floor function.
Thus we can apply Lemma \ref{lem:shape_1} to show that 
$\P \otimes \bP\hyphen \as$, the sum of the second and third terms 
of the right-hand side of \eqref{eq:shape_3} is smaller than 
\begin{align*}
 &d_\lambda ([\rho_1 nx],\lceil \rho_1 n/K \rceil Kx,\omega ,\sigma )
  +d_\lambda ([\rho_2 nx],\lfloor \rho_2 n/K \rfloor Kx,\omega ,\sigma )\\
 &< c_3(\lambda ) \epsilon (|[\rho_1 nx]|+|[\rho_2 nx]|).
\end{align*}
It follows that 
\begin{align*}
 &\limsup_{n \to \infty} \frac{1}{n} a_\lambda ([\rho_1 nx],[\rho_2 nx],\omega ,\sigma )\\
 &\leq \frac{\rho_2 -\rho_1}{K} \,\E \otimes \bE [a_\lambda (0,Kx,\omega ,\sigma )]
       +c_3(\lambda )\epsilon(\rho_1 +\rho_2 )|x|
\end{align*}
and therefore letting $\epsilon \searrow 0$ and $K \to \infty$, we obtain from Theorem \ref{thm:intro_1} that 
\begin{align}
 \limsup_{n \to \infty} \frac{1}{n} a_\lambda ([\rho_1 nx],[\rho_2 nx],\omega ,\sigma )
 &\leq (\rho_2 -\rho_1 )\alpha_\lambda (x)
 \label{eq:shape_4}
\end{align}
holds $\P \otimes \bP \hyphen \as$ 
Using Lemma \ref{lem:lyapunov_1}, we have 
\begin{align*}
 &a_\lambda ([\rho_1 nx],[\rho_2 nx],\omega ,\sigma )-a_\lambda ([\rho_1 x_n],[\rho_2 x_n],\omega ,\sigma )\\
 &\leq a_\lambda ([\rho_1 nx],[\rho_1 x_n],\omega ,\sigma )
       +a_\lambda ([\rho_2 x_n],[\rho_2 nx],\omega ,\sigma )\\
 &\leq d_\lambda ([\rho_1 nx],[\rho_1 x_n],\omega ,\sigma )
       +d_\lambda ([\rho_2 x_n],[\rho_2 nx],\omega ,\sigma ),
\end{align*}
and similarly
\begin{align*}
 &a_\lambda ([\rho_1 nx],[\rho_2 nx],\omega ,\sigma )-a_\lambda ([\rho_1 x_n],[\rho_2 x_n],\omega ,\sigma )\\
 &\geq -d_\lambda ([\rho_1 nx],[\rho_1 x_n],\omega ,\sigma )
       -d_\lambda ([\rho_2 nx],[\rho_2 x_n],\omega ,\sigma ).
\end{align*}
Furthermore, we have for $\epsilon \in \Q \cap (0,1/2)$ and sufficiently large $n$, 
\begin{align*}
 |[\rho_i nx]-[\rho_i x_n]| \leq 2\epsilon |[\rho_i nx]|,\qquad i=1,2.
\end{align*}
Lemma \ref{lem:shape_1} thereby implies that 
\begin{align*}
 &\limsup_{n \to \infty} \frac{1}{n}
  |a_\lambda ([\rho_1 nx],[\rho_2 nx],\omega ,\sigma )-a_\lambda ([\rho_1 x_n],[\rho_2 x_n],\omega ,\sigma )|\\
 &\leq \limsup_{n \to \infty}
       \frac{1}{n} \bigl( d_\lambda ([\rho_1 nx],[\rho_1 x_n],\omega ,\sigma )
                          +d_\lambda ([\rho_2 nx],[\rho_2 x_n],\omega ,\sigma ) \bigr) \\
 &\leq 2c_3(\lambda ) \epsilon (\rho_1 +\rho_2 )|x|
\end{align*}
holds $\P \otimes \bP \hyphen \as$ 
This together with \eqref{eq:shape_4} proves \eqref{eq:shape_2} since $\epsilon$ is arbitrary. 
\end{proof}
\begin{remark}
Zerner proved a stronger version of the shape theorem and used it to prove 
the large deviation lower bound in~\cite{Zer98b}. 
We find difficulty in proving a shape theorem in such a general 
form. Note that our Lyapunov exponent can be regarded as a mixture
of those in~\cite{Zer98a} and~\cite{Zer98b} and 
in the former paper, the {\em uniform 
shape theorem} requires some assumptions.
Our strategy, using the above directionally uniform shape theorem,
dates back to Sznitman's work on large deviations for Brownian
motion in Poissonian obstacles~\cite{Szn94b}. 
\end{remark}
\section{Large deviation estimates
}
Our goal in this section is to show Theorem \ref{thm:intro_3}. 
We prove upper and lower bounds of Theorem~\ref{thm:intro_3} 
in Subsections~\ref{LDP_upper} and~\ref{LDP_lower}, respectively.

\subsection{Upper bound}\label{LDP_upper}
In this subsection, we prove the upper bound \eqref{eq:intro_5} of Theorem~\ref{thm:intro_3}. 
Let us first mention some properties of the rate function $I$. 
We denote the essential domain of the rate function $I$ by
$\mathcal{D}_I$, that is, $\mathcal{D}_I:=\{ x \in \R^d;I(x)<\infty \}$. 
It is easy to see that $I$ is convex on $\R^d$, lower semicontinuous on
$\mathcal{D}_I$ and continuous on the interior 
of $\mathcal{D}_I$. 
Furthermore, by Theorem \ref{thm:intro_1}, we have an upper bound
\begin{equation*}
 I(x) \leq\sup_{\lambda\ge 0}\left[|x| \Bigl(\max_{|e|=1} 
\E[-\log\omega (0,e)]+
\bE[\theta_{\lambda,\sigma}(0)] \Bigr)-\lambda\right].
\end{equation*}
Since 
$\bE[\theta_{\lambda,\sigma}(0)]\sim \lambda\bE[\int_0^\infty s\,\sigma_0(ds)]$ 
as $\lambda\downarrow 0$ and is concave in $\lambda$, the above supremum
is attained at $\lambda=0$ if 
$|x| \leq \bE[\int_0^\infty s\,\sigma_0(ds)]^{-1}$ and whence
\begin{align}
 0 \leq I(x) \leq |x| \max_{|e|=1} \E[-\log \omega (0,e)].
 \label{eq:ldpupper_3}
\end{align}
\begin{proof}[\bf Proof of the upper bound in Theorem \ref{thm:intro_3}]
Let us first show that 
\begin{align}
 \lim_{R \to \infty} \limsup_{t \to \infty}
 \frac{1}{t} \log \tilde{P}_{\omega,\sigma}^0 (Z_t \not\in t\bar{B}(0,R))=-\infty
 \label{eq:ldpupper_1}
\end{align}
holds $\P \otimes \bP \hyphen \as$ 
We have for any $\lambda,t \geq 0$ and subset $K \subset \R^d\setminus\{0\}$, 
\begin{align*}
 \tilde{P}_{\omega,\sigma}^0 (Z_t \in tK)
 &\leq \exp \{ \lambda t \}
       \tilde{E}_{\omega ,\sigma}^0 \bigl[ \exp \{ -\lambda t \} \1{\{ H^Z([tK]) \leq t \}} \bigr] \\
 &\leq \exp \{ \lambda t \} \,e_\lambda (0,[tK],\omega ,\sigma ),
\end{align*}
and hence Corollary \ref{cor:shape_1} implies that 
\begin{align}
\begin{split}
 \limsup_{t \to \infty} \frac{1}{t} \log \tilde{P}_{\omega ,\sigma}^0(Z_t \in tK)
 &\leq \limsup_{t \to \infty} \biggl( \lambda -\frac{1}{t} \log e_\lambda (0,[tK],\omega ,\sigma ) \biggr) \\
 &= \lambda -\inf_{x \in K} \alpha_\lambda (x)
\end{split}
\label{eq:ldpupper_2}
\end{align}
holds $\P \otimes \bP \hyphen \as$ 
By taking $K=\bar{B}(0,R)^c$, it follows from Theorem \ref{thm:intro_1} that we obtain $\P \otimes \bP \hyphen \as$, 
\begin{align*}
 \limsup_{t \to \infty} \frac{1}{t} \log \tilde{P}_{\omega ,\sigma}^0 (Z_t \not\in t\bar{B}(0,R))
 &\leq -\inf_{x \in \bar{B}(0,R)^c}(\alpha_\lambda (x)-\lambda )\\
 &\leq -\inf_{x \in \bar{B}(0,R)^c}(|x|(-\log \bE[\exp \{ -\theta_{\lambda ,\sigma}(0) \} ])-\lambda )\\
 &\leq -R (-\log \bE[\exp \{ -\theta_{\lambda ,\sigma}(0) \} ])+\lambda ,
\end{align*}
which proves \eqref{eq:ldpupper_1} by letting $R \to \infty$. 

Now we show the upper bound \eqref{eq:intro_5}. 
It suffices to consider compact $A\subset\R^d$ thanks to~\eqref{eq:ldpupper_1}. 
Moreover, we may assume $0 \not\in A$ since $\inf_{x \in A} I(x)=0$ if $0 \in A$ by \eqref{eq:ldpupper_3}. 
For every $\delta >0$ we introduce the $\delta$-rate function $I^\delta$ as
\begin{align*}
 I^\delta (x):=(I(x)-\delta ) \wedge \frac{1}{\delta}
\end{align*}
and set  
\begin{align*}
 A_\lambda (\delta ):=\Bigl\{ y \in A;\alpha_\lambda (y)-\lambda >\inf_{x \in A}I^\delta (x)\Bigr\} 
\end{align*}
for each $\lambda \geq 0$. 
Applying \eqref{eq:ldpupper_2} with $K=A_\lambda (\delta)$, we obtain from Corollary \ref{cor:shape_1} that 
\begin{align*}
 \limsup_{t \to \infty} \frac{1}{t} \log \tilde{P}_{\omega ,\sigma}^0(Z_t \in tA_\lambda (\delta ))
 \leq \lambda -\inf_{x \in A_\lambda (\delta )} \alpha_\lambda (x)
 \leq -\inf_{x \in A} I^\delta (x).
\end{align*}
Since $A=\bigcup_{\lambda \geq 0}A_\lambda (\delta)$ and 
$A$ is compact, there are $\lambda_i$ ($1 \leq i \leq m$) such that $A_{\lambda_i}(\delta)$ ($1 \leq i \leq m$) cover $A$. 
Thus, for any $\delta >0$, 
\begin{align*}
 \limsup_{t \to \infty} \frac{1}{t} \log \tilde{P}_{\omega ,\sigma}^0(Z_t \in tA)
 &\leq \max_{1 \leq i \leq m}
      \limsup_{t \to \infty} \frac{1}{t} \log \tilde{P}_{\omega ,\sigma}^0(Z_t \in tA_{\lambda_i}(\delta ))\\
 &\leq  -\inf_{x \in A}I^\delta (x).
\end{align*}
Since $\lim_{\delta \searrow 0} \inf_{x \in A} I^\delta (x)=\inf_{x \in A}I(x)$, the upper bound \eqref{eq:intro_5} follows by letting $\delta \searrow 0$. 
\end{proof}
\subsection{Lower bound}\label{LDP_lower}
In this subsection, we prove the lower bound \eqref{eq:intro_6} of Theorem~\ref{thm:intro_3}. 
Let us start with the following lemma. 

\begin{lemma}\label{lem:ldplower_0}
Let $x \in \Q^d \setminus \{ 0 \}$ and assume that $\rho_1,\rho_2 \in \R$ satisfy $0 \leq \rho_1<\rho_2$. 
Then the following holds $\P \otimes \bP \hyphen \as$:
\begin{align}
 \lim_{t \to \infty} \hat{P}_{\lambda ,\omega ,\sigma}^{[\rho_1 x_t],[\rho_2 x_t]}
 \biggl( \frac{H^Z([\rho_2 x_t])}{(\rho_2 -\rho_1 )t} \in (\gamma_1 ,\gamma_2 ) \biggr) =1
 \label{eq:ldplower_1}
\end{align}
for all sequences $(x_t)_{t \geq 0}$ of $\R^d$ with $x_t/t \to x$ as $t \to \infty$ and all $\lambda >0,\,\gamma_1,\gamma_2 \in \R$ satisfying 
\begin{align}
 0 \leq \gamma_1 <\alpha'_{\lambda +}(x) \leq \alpha'_{\lambda -}(x)<\gamma_2 .
 \label{eq:ldplower_2}
\end{align}
\end{lemma}
\begin{proof}
Corollary \ref{cor:shape_2} implies that $\P \otimes \bP \hyphen \as$, for $x,\rho_1,\rho_2,(x_t)_{t \geq 0},\lambda,\gamma_1,\gamma_2$ as above and $\mu \in (0,\lambda)$, 
\begin{align*}
 &\limsup_{t \to \infty} \frac{1}{t} \log \hat{P}^{[\rho_1 x_t],[\rho_2 x_t]}_{\lambda ,\omega ,\sigma}
                                          (H^Z([\rho_2 x_t]) \geq (\rho_2 -\rho_1 )\gamma_2 t)\\
 &= (\rho_2 -\rho_1 )\alpha_\lambda (x)\\
 &\quad
    +\limsup_{t \to \infty} \frac{1}{t}
     \log \tilde{E}_{\omega ,\sigma}^{[\rho_1 x_t]} \bigl[ \exp \{ (-\lambda +\mu )H^Z([\rho_2 x_t]) \}
                                         \exp \{ -\mu H^Z([\rho_2 x_t]) \}\\
 &\qquad \qquad \qquad \qquad \qquad \ \ \times
     \1{\{ H^Z([\rho_2 x_t])<\infty ,\,H^Z([\rho_2 x_t]) \geq (\rho_2 -\rho_1 )\gamma_2 t \}} \bigr] \\
 &\leq (\rho_2 -\rho_1 )\alpha_\lambda (x)-\mu (\rho_2 -\rho_1 )\gamma_2 -(\rho_2 -\rho_1 )\alpha_{\lambda -\mu}(x)\\
 &= \mu (\rho_2 -\rho_1 ) \biggl( \frac{\alpha_\lambda (x)-\alpha_{\lambda -\mu}(x)}{\mu} -\gamma_2 \biggr) .
\end{align*}
It follows from \eqref{eq:ldplower_2} that the most right-hand side of the above expression is negative for $\mu$ small enough. 
Since the corresponding statement can be proved for the event $\{
 H^Z([\rho_2 x_t]) \leq (\rho_2 -\rho_1 )\gamma_1 t \}$ in the same manner, 
we have \eqref{eq:ldplower_1}. 
\end{proof}

\begin{proof}[\bf Proof of the lower bound in Theorem \ref{thm:intro_3}]
It suffices to show that $\P \otimes \bP \hyphen \as$, 
\begin{align}
 \liminf_{t \to \infty} \frac{1}{t} \log \tilde{P}_{\omega ,\sigma}^0(Z_t \in tB(z,r)) \geq -I(z)
 \label{eq:ldplower_3}
\end{align}
holds for all $z \in \Q^d \setminus \{ 0 \} \cap \mathcal{D}_I$ and $0<r \in \Q$. 
To prove this, let us define 
\begin{align*}
 \lambda_\infty:=
  \sup \{ \lambda >0;\alpha'_{\lambda -}(z) \geq 1 \},
\end{align*}
with the convention $\sup \emptyset=0$. 
It is easy to check that $I(z)=\alpha_{\lambda_\infty}(z)-\lambda_\infty$ in the case $\lambda_\infty<\infty$, and $I(z)=\lim_{\lambda \to \infty}(\alpha_\lambda (z)-\lambda)$ otherwise. 
We first treat the case $\lambda_\infty<\infty$. 
By the concavity of $\alpha_{\lambda}(z)$ in $\lambda$, we can find
 sequences $(\gamma_n)_{n=1}^\infty,(\delta_n)_{n=1}^\infty$ and
 $(\lambda_n)_{n=1}^\infty$ such that 
\begin{itemize}
 \item if $\alpha'_{\lambda -}(z)<1$ for all $\lambda >0$, then $\alpha'_{\lambda_n}(z)$ exists, $\lambda_n \to 0$ and 
       \begin{align*}
        \gamma_n :=\alpha'_{\lambda_n}(z) \Bigl( 1-\frac{1}{n} \Bigr) ,\quad
        \delta_n :=\alpha'_{\lambda_n}(z) \biggl( 1+\frac{1-\alpha'_{\lambda_n}(z)}{n} \biggr) <1,
       \end{align*}
\item otherwise, $\lambda_n \to \lambda_\infty$, $1-2/n \in [\alpha'_{\lambda_n+}(z),\alpha'_{\lambda_n-}(z)]$ and 
       \begin{align*}
        \gamma_n :=1-\frac{3}{n},\quad \delta_n :=1-\frac{1}{n}.
       \end{align*}
\end{itemize}
Observe that for the above sequences, we have
       \begin{align}
        (\gamma_n ,\delta_n ) \cap [\alpha'_{\lambda_n+}(z),\alpha'_{\lambda_n-}(z)] \not= \emptyset .
        \label{eq:ldplower_4}
       \end{align}
Now recall that Assumption~\ref{assumption}-(2) is equivalent to the 
following (see~\cite[Proposition 8]{Zer98b}): 
for all $\epsilon >0$ there is some $R(\epsilon) \geq 2$ such that 
\begin{align}
 \P \bigl( P_\omega^0 (X_{R(\epsilon )}=0)>e^{-\epsilon R(\epsilon )} \bigr) >0.
 \label{eq:intro_7}
\end{align}
We choose $R(\epsilon)>0$ satisfying \eqref{eq:intro_7} for $\epsilon >0$ and fix $\delta >0$ with $\bP(\theta_{1,\sigma}(0) \geq \delta)>0$. 
Then, for $y \in (2R(\epsilon)+1)\Z^d$
\begin{align*}
\Phi_0 (y,\epsilon )
:= \Bigl\{(\omega,\sigma): 
P_\omega^y (X_{R(\epsilon )}=y)>e^{-\epsilon R(\epsilon )},
\min_{x \in B(y,R(\epsilon ))} \theta_{1,\sigma}(x) \geq \delta \Bigr\}
\end{align*}
has strictly positive $\P \otimes \bP$-probability and is independent for
distinct $y$'s. 
Let $y_t=y_t(\epsilon, \omega ,\sigma ) \in
 (2R(\epsilon)+1)\Z^d$ be a vertex with minimal distance from 
$[tz]$ 
such that $(\omega,\sigma) \in \Phi_0(y_t,\epsilon)$. 
A simple application of the Borel--Cantelli lemma shows that $\P \otimes \bP \hyphen \as$, these $y_t$ exist and satisfy 
\begin{align}
 |tz-y_t| \leq 2(\log t)^2
 \label{eq:ldplower_5}
\end{align}
for all sufficiently large $t$. 
Let us introduce additional notation. 
Denote 
\begin{align*}
 T_n^X(\epsilon ,t):=\inf \Biggl\{ k \geq 0;\sum_{m=0}^{k-1} \tau_m (X_m)>\gamma_n t,\,X_k=y_t \Biggr\}
\end{align*}
and 
\begin{align*}
 T_n^Z(\epsilon ,t):=\inf \{ s>\gamma_n t;Z_s=y_t \} = \sum_{m=0}^{T_n^X(t)-1} \tau_m (X_m).
\end{align*}
We then define the random variable
\begin{align*}
 b_n(\epsilon ,t,\omega ,\sigma )
 := -\log \tilde{E}_{\omega ,\sigma}^0 \bigl[ \exp \{ -\lambda_n T_n^Z(\epsilon ,t) \}
    \1{\{ T_n^Z(\epsilon ,t)<\delta_n t \}} \bigr]
\end{align*}
and the event
\begin{align*}
 \Lambda_1 (n,\epsilon ,t):= 
\{(X, \tau): Z_{s+T_n^Z(\epsilon ,t)} \in B(y_t,R(\epsilon)) \textrm{ for all $s \in [0,(1-\gamma_n)t]$} \} .
\end{align*}

Now, the left-hand side of \eqref{eq:ldplower_3} is greater than 
\begin{align*}
 \lambda_n \gamma_n +\liminf_{t \to \infty} \frac{1}{t}
 \log \tilde{E}_{\omega ,\sigma}^0 \bigl[ \exp \{ -\lambda_n T_n^Z(\epsilon ,t) \}
 \1{\{ T_n^Z(\epsilon ,t)<\delta_n t \} \cap \Lambda_1 (n,\epsilon ,t)} \bigr]
\end{align*}
since $B(y_t, R(\epsilon))\subset tB(z,r)$ for sufficiently large
 $t$ by \eqref{eq:ldplower_5}. 
The strong Markov property shows that the above expression equals to 
\begin{align}
 \lambda_n \gamma_n
 -\limsup_{t \to \infty} \frac{1}{t} b_n(\epsilon,t, \omega ,\sigma )
 +\liminf_{t \to \infty} \frac{1}{t}
  \log \sum_{\ell =0}^\infty \tilde{P}_{\omega ,\sigma}^{y_t}(\Lambda_2 (n,\epsilon ,t,\ell )),
 \label{eq:ldplower_6}
\end{align}
where $\Lambda_2 (n,\epsilon ,t,\ell )$ is the event defined as 
\begin{align*}
 \Lambda_2 (n,\epsilon ,t,\ell )
 := \Biggl\{ (X, \tau):
&\sum_{m=0}^{\ell -1} \tau_m (X_m) \leq (1-\gamma_n )t<\sum_{m=0}^\ell \tau_m (X_m),\\
&X_m \in B(y_t,R(\epsilon )) \textrm{ for all $m \in [0,\ell-1]$} \Biggr\}.
\end{align*}
To control the second and third terms of \eqref{eq:ldplower_6}, we use the following two lemmas. 
\begin{lemma} \label{lem:ldplower_2}
For any $\epsilon >0$ and $n \geq 1$, we have $\P \otimes \bP \hyphen \as$, 
\begin{align*}
 \lim_{t \to \infty} \frac{1}{t} b_n(\epsilon ,t,\omega ,\sigma )=\alpha_{\lambda_n}(z).
\end{align*}
\end{lemma}

\begin{lemma} \label{lem:ldplower_3}
For any $\epsilon >0$ and $n \geq 1$, we have $\P \otimes \bP \hyphen \as$, 
\begin{align}
 \liminf_{t \to \infty} \frac{1}{t}
  \log \sum_{\ell =0}^\infty \tilde{P}_{\omega ,\sigma}^{y_t}(\Lambda_2 (n,\epsilon ,t,\ell ))
 \geq -\frac{2\epsilon }{\delta}.
 \label{eq:ldplower_10}
\end{align}
\end{lemma}
Let us postpone the proofs of these lemmas to the end of this
 subsection. 
It follows from Lemmas \ref{lem:ldplower_2}, \ref{lem:ldplower_3}, and
from \eqref{eq:ldplower_6} that $\P \otimes \bP \hyphen \as$, 
\begin{align*}
 \liminf_{t \to \infty} \frac{1}{t} \log \tilde{P}_{\omega ,\sigma}^0(Z_t \in tB(z,r))
 \geq \lambda_n \gamma_n -\alpha_{\lambda_n}(z)-\frac{2\epsilon }{\delta},
\end{align*}
which completes the proof of \eqref{eq:ldplower_3} in the case $\lambda_\infty <\infty$ by letting $\epsilon \searrow 0$ and $n \to \infty$. 

We next treat the case $\lambda_\infty=\infty$. 
In this case, $\alpha'_{\lambda -}(uz)=u \alpha'_{\lambda -}(z)<1$ holds for all $u \in \Q \cap (0,1)$ and all sufficiently large $\lambda$. 
Moreover, for $u \in (0 \vee (1-r/|z|),1)$ we pick $0<r'(u) \in \Q$ with $B(u z,r'(u)) \subset B(z,r)$. 
Applying the same argument as in the case $\lambda_\infty <\infty$ and
 using the convexity of the rate function $I$, one can show that 
\begin{align*}
 &\liminf_{t \to \infty} \frac{1}{t} \log \tilde{P}_{\omega ,\sigma}^0(Z_t \in tB(z,r))\\
 &\geq \liminf_{t \to \infty} \frac{1}{t} \log \tilde{P}_{\omega ,\sigma}^0(Z_t \in tB(u z,r'(u)))
  \geq -I(u z) \geq -u I(z)
\end{align*}
holds $\P \otimes \bP \hyphen \as$ 
This proves \eqref{eq:ldplower_3} by letting $u \nearrow 1$. 
\end{proof}

We close this section with the proofs of Lemmas \ref{lem:ldplower_2} 
and \ref{lem:ldplower_3}.
\begin{proof}[\bf Proof of Lemma \ref{lem:ldplower_2}]
Note that $b_n(\epsilon,t,\omega,\sigma) \geq
 a_{\lambda_n}(0,y_t,\omega,\sigma)$ since $H^Z(y_t) \leq T_n^Z(\epsilon,t)$. 
Theorem \ref{thm:intro_2} hence implies that we have $\P \otimes \bP \hyphen \as$, 
\begin{align*}
 \alpha_{\lambda_n}(z) \leq \liminf_{t \to \infty} \frac{1}{t} b_n(\epsilon ,t,\omega ,\sigma ).
\end{align*}
It remains to show that 
\begin{align}
 \limsup_{t \to \infty} \frac{1}{t} b_n(\epsilon ,t,\omega ,\sigma ) \leq \alpha_{\lambda_n}(z)
 \label{eq:ldplower_9}
\end{align}
holds $\P \otimes \bP \hyphen \as$ 
Thanks to \eqref{eq:ldplower_4}, we can pick $\rho \in (0,1)$ and $\eta >0$ such that 
\begin{align*}
 \rho \alpha'_{\lambda_n +}(z)+(1-\rho )\alpha'_{\lambda_n -}(z)+[-\eta ,\eta ] \subset (\gamma_n ,\delta_n ).
\end{align*}
Setting $\xi_t:=[\rho y_t]$, we know from the choice of $\rho$ and $\eta$ that $T_n^Z(\epsilon,t)<\delta_n t$ holds on 
\begin{align*}
 \Lambda_3 (n,\epsilon ,t)
 :=\Biggl\{ (X,\tau): 
&\frac{1}{\rho t} \sum_{m=0}^{H^X(\xi_t)-1} \tau_m (X_m) \in \alpha'_{\lambda_n +}(z)+[-\eta ,\eta ],\\
            &\frac{1}{(1-\rho )t} \sum_{m=H^X(\xi_t)}^{H_{\xi_t}^X(y_t)-1} \tau_m (X_m)
             \in \alpha'_{\lambda_n -}(z)+[-\eta ,\eta ] \Biggr\} .
\end{align*}
It follows from this and the strong Markov property that 
\begin{align*}
 &\tilde{E}_{\omega ,\sigma}^0 \bigl[ \exp \{ -\lambda_n T_n^Z(\epsilon ,t) \}
  \1{\{ T_n^Z(\epsilon ,t)<\delta_n t \}} \bigr] \\
 &\geq \tilde{E}_{\omega ,\sigma}^0
       \Biggl[ \exp \Biggl\{ -\lambda_n \sum_{m=0}^{H^X_{\xi_t}(y_t)-1} \tau_m (X_m) \Biggr\}
               \1{\Lambda_3 (n,\epsilon ,t)} \Biggr] \\
 &= \tilde{E}_{\omega ,\sigma}^0 \bigl[ \exp \{ -\lambda_n H^Z(\xi_t) \}
    \1{\{ H^Z(\xi_t)/(\rho t) \in \alpha'_{\lambda_n +}(z)+[-\eta ,\eta ] \}} \bigr] \\
 &\quad \times
    \tilde{E}_{\omega ,\sigma}^{\xi_t} \bigl[ \exp \{ -\lambda_n H^Z(y_t) \}
    \1{\{ H^Z(y_t)/((1-\rho )t) \in \alpha'_{\lambda_n -}(z)+[-\eta ,\eta ] \}} \bigr].
\end{align*}
Let $\mu_1$ and $\mu_2$ be such that $0<\mu_1<\lambda_n <\mu_2$
and $\alpha'_{\lambda_n +}(z)-\eta<\alpha'_{\mu_1 +}(z)<\alpha'_{\mu_2 -}(z)<\alpha'_{\lambda_n -}(z)+\eta$. Then the most right-hand side of the above expression is greater than 
\begin{align*}
 &e_{\mu_2}(0,\xi_t,\omega ,\sigma )
  \,\hat{P}_{\mu_2 ,\omega ,\sigma}^{0,\xi_t}
  \biggl( \frac{H^Z(\xi_t)}{\rho t} \in \alpha'_{\lambda_n +}(z)+[-\eta ,\eta ] \biggr) \\
 &\times \exp \{ -(\lambda_n -\mu_1 )(\alpha'_{\lambda_n -}(z)+\eta )(1-\rho )t \} \\
 &\times e_{\mu_1}(\xi_t ,y_t,\omega ,\sigma )
         \,\hat{P}_{\mu_1 ,\omega ,\sigma}^{\xi_t ,y_t}
         \biggl( \frac{H^Z(y_t)}{(1-\rho )t} \in \alpha'_{\lambda -}(x)+[-\eta ,\eta ] \biggr) .
\end{align*}
We thereby obtain for $t>0$, 
\begin{align*}
 \frac{1}{t} b_n(\epsilon ,t,\omega ,\sigma )
 &\leq \frac{1}{t} a_{\lambda_2}(0,\xi_t,\omega ,\sigma )
       -\frac{1}{t} \log \hat{P}_{\mu_2 ,\omega ,\sigma}^{0,\xi_t}
        \biggl( \frac{H^Z(\xi_t)}{\rho t} \in \alpha'_{\lambda_n +}(z)+[-\eta ,\eta ] \biggr)\\
 &\quad
       +(\lambda_n -\mu_1 )(\alpha'_{\lambda_n -}(z)+\eta )
       +\frac{1}{t} a_{\lambda_1}(\xi_t ,y_t,\omega ,\sigma )\\
 &\quad
       -\frac{1}{t} \log \hat{P}_{\lambda_1 ,\omega ,\sigma}^{\xi_t ,y_t}
        \biggl( \frac{H^Z(y_t)}{(1-\rho )t} \in \alpha'_{\lambda -}(z)+[-\eta ,\eta ] \biggr) .
\end{align*}
Note that we have $y_t/t \to z$ from \eqref{eq:ldplower_5}. 
Therefore, applying Corollary~\ref{cor:shape_2} and Lemma~\ref{lem:ldplower_0}, we get $\P \otimes \bP \hyphen \as$, 
\begin{align*}
 \limsup_{t \to \infty} \frac{1}{t} b_n(\epsilon ,t,\omega ,\sigma )
 \leq \rho \alpha_{\mu_2}(z)+(\lambda_n -\mu_1 )(\alpha'_{\lambda_n -}(z)+\eta )
      +(1-\rho ) \alpha_{\mu_1}(z),
\end{align*}
which concludes \eqref{eq:ldplower_9} by letting $\mu_1 \nearrow \lambda_n$ and $\mu_2 \searrow \lambda_n$. 
\end{proof}

\begin{proof}[\bf Proof of Lemma \ref{lem:ldplower_3}]
Let $L(n,t):=\lfloor 2t/\delta \rfloor +1$. 
If $\sum_{m=0}^{L(n,t)-1} \tau_m (X_m) \geq \delta L(n,t)/2$, then $\sum_{m=0}^{L(n,t)-1} \tau_m (X_m)>(1-\gamma_n)t$. 
Thus it follows that 
\begin{align*}
 &\sum_{\ell =0}^\infty \tilde{P}_{\omega ,\sigma}^{y_t}(\Lambda_2 (n,\epsilon ,t,\ell ))\\
 &\geq \sum_{\ell =0}^{L(n,t)-1} \tilde{P}_{\omega ,\sigma}^{y_t} \Biggl(
       \sum_{m=0}^{\ell -1} \tau_m (X_m) \leq (1-\gamma_n )t<\sum_{m=0}^\ell \tau_m (X_m),\\
 &\qquad \qquad \qquad \quad
       \sum_{m=0}^{L(n,t)-1} \tau_m (X_m) \geq \frac{\delta}{2}L(n,t),\\
 &\qquad \qquad \qquad \qquad
       X_m \in B(y_t,R(\epsilon )) \textrm{ for all $m \in [0,L(n,t)-1]$} \Biggr) \\
 &= \tilde{P}_{\omega ,\sigma}^{y_t} \Biggl(
    \sum_{m=0}^{L(n,t)-1} \tau_m (X_m) \geq \frac{\delta}{2}L(n,t),\,
    X_k \in B(y_t,R(\epsilon )) \textrm{ for all $k \in [0,L(n,t)-1]$} \Biggr) \\
 &= E_\omega^{y_t} \Biggl[ P_\sigma^\HT \Biggl( \sum_{m=0}^{L(n,t)-1} \tau_m (X_m) \geq \frac{\delta}{2}L(n,t) \Biggr) 
    \1{\{ X_k \in B(y_t,R(\epsilon )) \textrm{ for all $k \in [0,L(n,t)-1]$} \}} \Biggr] .
\end{align*}
On the other hand, the choice of $y_t$ and Chebyshev's inequality 
imply that 
\begin{align*}
 P_\sigma^\HT \Biggl( \sum_{m=0}^{L(n,t)-1} \tau_m (X_m) \geq \frac{\delta}{2}L(n,t) \Biggr)
 &\geq 1-e^{(\delta /2)L(n,t)}
       E_\sigma^\HT \Biggl[ \exp \Biggl\{ -\sum_{m=0}^{L(n,t)-1} \tau_m (X_m) \Biggr\} \Biggr]\\
 &\geq 1-e^{(\delta /2)L(n,t)}
       \prod_{m=0}^{L(n,t)-1} \exp \{ -\theta_{1,\sigma}(X_m) \}\\
 &\geq 1-e^{-(\delta /2)L(n,t)}
\end{align*}
uniformly in paths $(X_n)_{n=0}^\infty$ with $X_m \in B(y_t,R(\epsilon))$ for all $m \in [0,L(n,t)-1]$. 
It follows from the choice of $y_t$ (recall \eqref{eq:intro_7}) that the left-hand side of \eqref{eq:ldplower_10} is greater than 
\begin{align*}
 &\liminf_{t \to \infty} \frac{1}{t}
  \Bigl( \log \bigl( 1-e^{-(\delta /2)L(n,t)} \bigr)
         +\log P_\omega^{y_t}(X_{R(\epsilon )}=y_t)^{L(n,t)/R(\epsilon )} \Bigr) \\
 &\geq \liminf_{t \to \infty}
       \biggl( \frac{1}{t}\log \bigl( 1-e^{-(\delta /2)L(n,t)} \bigr)
               -\frac{2\epsilon }{\delta} \biggr) \\
 &= -\frac{2\epsilon }{\delta}.
\end{align*}
Since $\delta$ is fixed and $\epsilon$ is arbitrary,
this proves the lemma. 
\end{proof}
\section{First passage percolation}
In this section, we relate our Lyapunov exponent to the so-called 
{\it time constant} of a first passage percolation in the limit
$\lambda \to \infty$. 
This will be used in the next section to study the asymptotics of the 
rate function. 
Throughout this section, we assume that for some deterministic 
function $L(\lambda)$ with $\lim_{\lambda\to\infty}L(\lambda)=\infty$,
\begin{equation}
 \lim_{\lambda\to\infty}
 \frac{\theta_{\lambda,\sigma}(z)}{L(\lambda)}
 =\Theta_{\sigma}(z)\in (0,\infty) 
\textrm{ exists }\bP\textrm{-}\as
\label{eq:percolation_0}
\end{equation}
As the following examples show, this is a rather restrictive assumption.
\begin{ex}\label{ex:percolation_0}
Let us denote the distribution function of $\sigma(0)$ by $F_{\sigma}$.
\begin{enumerate}
 \item If $\inf \supp \sigma(0)>0$ for each $\sigma$, 
then~\eqref{eq:percolation_0} holds with $L(\lambda)=\lambda$. 
 \item If for each $\sigma$ there is a $\gamma(\sigma)>0$ such that 
$\lim_{x\downarrow 0}F_{\sigma}(x)/x^{\gamma(\sigma)} \in (0,\infty)$, 
then \eqref{eq:percolation_0} holds
with $L(\lambda)=\log \lambda$. This includes the case where
all $\{\sigma_x\}_{x\in\Z^d}$ are exponential distributions,
which is sometimes called the ``random hopping time dynamics''. 
 \item If $x^{\gamma(\sigma)}$ in the previous example is replaced
by $x^{\gamma(\sigma)}(\log x)^{\delta(\sigma)}$ with some 
non-constant $\delta(\sigma)$, then~\eqref{eq:percolation_0} 
fails to hold.
 \item If there exists a $\gamma>0$ such that 
$-\lim_{x\downarrow 0}\log F_{\sigma}(x)/x^{-\gamma}\in (0,\infty)$ 
for each $\sigma$, then~\eqref{eq:percolation_0} holds with 
$L(\lambda)=\lambda^{\gamma/(\gamma+1)}$.
 \item If there exist $\sigma_1$ and $\sigma_2$ such that 
$-\lim_{x\downarrow 0}\log F_{\sigma_i}(x)/x^{-\gamma(\sigma_i)}\in (0,\infty)$ 
for $i=1,2$ with different positive constants $\gamma(\sigma_1)$
and $\gamma(\sigma_2)$, 
then~\eqref{eq:percolation_0} fails to hold.
\end{enumerate}
These are well-known facts in Tauberian theory, see~\cite{BGT89}. 
\hfill$\square$
\end{ex}
 
For given positive i.i.d.~random variables $\{\xi(z)\}_{z\in\Z^d}$, 
we define the passage time of a nearest neighbor path
$r=(r_0,r_1,\dots,r_n)$ as 
\begin{align*}
 T(r,\xi):=\sum_{i=0}^{n-1} \xi(r_i),
\end{align*}
where the right hand side is set to be 0 if $n=0$. 
The travel time from $x$ to $y$ is defined as 
\begin{align*}
 T(x,y,\xi):=\inf \{ T(r,\xi);
\textrm{$r$ is a path from $x$ to $y$} \}.
\end{align*}
It is shown by Cox and Durrett~\cite{CD81}
that there exists a deterministic norm $\nu_{\xi}$ such that 
\begin{align}
\begin{split}
 \lim_{n \to \infty} \frac{1}{n} \,T(0,nx,\xi)
  = \nu_{\xi} (x) \textrm{ in probability} 
\end{split}
\label{eq:percolation_1}
\end{align}
for all $x \in \Z^d$. 
\begin{proposition} \label{prop:percolation_2}
For any $x \in \R^d$, 
\begin{align}
 \frac{\alpha_\lambda (x)}{L(\lambda)} \to \nu_{\Theta_\sigma} (x),
 \quad \lambda \to \infty .
 \label{eq:percolation_4}
\end{align}
\end{proposition}
\begin{proof}
It suffices show the assertion only for $x\in\Z^d$. 
Indeed, since both $\alpha_{\lambda}$ and $\nu_{\Theta_\sigma}$ are 
homogeneous, 
it extends to $\Q^d$ and then to $\R^d$ by continuity. 
Now for $x\in\Z^d$, we estimate the difference as 
\begin{equation}
\begin{split}
\left|\frac{\alpha_\lambda (x)}{L(\lambda)}-\nu_{\Theta_\sigma} (x)\right|
&\le \left|\frac{\alpha_\lambda (x)}{L(\lambda)}-
\frac{a_\lambda (0,nx,\omega ,\sigma )}{nL(\lambda)}\right|\\
&\quad+\left|\frac{a_\lambda (0,nx,\omega ,\sigma )}{nL(\lambda)}
-\frac{1}{n}T(0,nx,\theta_{\lambda,\sigma}/L(\lambda))\right|\\
&\quad+\left|\frac{1}{n}T(0,nx,\theta_{\lambda,\sigma}/L(\lambda))
-\nu_{\theta_{\lambda,\sigma}/L(\lambda)}(x)\right|\\
&\quad+\left|\nu_{\theta_{\lambda,\sigma}/L(\lambda)}(x)-\nu_{\Theta_\sigma}(x)\right|,
\label{eq: percolation_4-1}
\end{split}
\end{equation}
where $n\in\N$. 
Note that for any fixed $\lambda>0$, 
the first and third terms converge to 0 in probability as 
$n\to\infty$. 
We also know that the fourth term in~\eqref{eq: percolation_4-1}
tends to $0$ as $\lambda\to\infty$ due to our 
assumption~\eqref{eq:percolation_0} and the continuity 
of the time constant shown in Theorem~6.9 in~\cite{Kes86}. 

The following lemma gives a control on the second term. 
\begin{lemma}\label{lem:percolation_3}
For any $\epsilon>0$, 
there exists a $\Lambda>0$ such that for all $\lambda\ge \Lambda$, 
\begin{equation}
\limsup_{n\to\infty}\P\otimes\bP
\left(\left|\frac{a_\lambda (0,nx,\omega ,\sigma )}{nL(\lambda)}
-\frac{1}{n}T(0,nx,\theta_{\lambda,\sigma}/L(\lambda))\right|>\epsilon
\right)<\epsilon.
\label{uniform}
\end{equation} 
\end{lemma}
\begin{proof}
Since one of the bounds, namely
\begin{align*}
 \frac{a_\lambda (0,nx,\omega ,\sigma )}{L(\lambda)}
 &\geq -\frac{1}{L(\lambda)} \log E^0_\omega 
 [\exp \{ - T(0,nx,\theta_{\lambda, \sigma}) \} 
\,\1{\{ H^X(nx)<\infty \}}]\\
 &\geq T(0,nx,\theta_{\lambda,\sigma}/L(\lambda))
\end{align*}
is trivial, we have only to show that for any $\epsilon >0$, 
\begin{align*}
\P\otimes\bP
\left(\frac{a_\lambda (0,nx,\omega ,\sigma)}
{nL(\lambda)} \leq \frac{1}{n}
T(0,nx,\theta_{\lambda,\sigma}/L(\lambda))+\epsilon\right)
>1-\epsilon 
\end{align*}
when $\lambda$ and $n$ are sufficiently large.
To this end, we first pick a path $r=\{r_m\}_{m=0}^{N(r)}$ from 
those paths connecting $0$ and $nx$ and satisfying  
\begin{equation}
 T(r,\theta_{\lambda, \sigma}/L(\lambda))
 \le T(0,nx,\theta_{\lambda, \sigma}/L(\lambda))+1
\label{near-min}
\end{equation}
by some deterministic rule. 
\begin{lemma}\label{percolation:lem_4}
Fix $x \in \Z^d$. For sufficiently large $\lambda$, 
there exists a constant $c_x>0$ such that 
\begin{equation}
 \lim_{n\to\infty}\bP\left(N(r)\le c_x n\right)=1,
\end{equation}
where $N(r)$ is the length of the path $r$ picked above. 
\end{lemma}
\begin{proof}
Note first that for sufficiently large $\lambda$, 
\begin{equation*}
\lim_{n\to\infty}
\bP(T(0,nx,\theta_{\lambda, \sigma}/L(\lambda))\le 2\nu_{\Theta_\sigma}(x)n)
=1
\end{equation*}
by~\eqref{eq:percolation_1} and the continuity of the time constant. 
Hence 
\begin{equation}
\begin{split}
 \limsup_{n\to\infty}\bP\left(N(r)>cn\right)
 &\le \limsup_{n\to\infty}\bP\left(T(r,\theta_{\lambda, \sigma}/L(\lambda))
 \le2\nu_{\Theta_\sigma}(x)n\right)\\
 &\le \limsup_{n\to\infty}\exp\left\{2\nu_{\Theta_\sigma}(x)n+\sum_{m=1}^{cn}
 \log\bE[e^{-\theta_{\lambda,\sigma}(r_m)/L(\lambda)}]\right\}.
\end{split}
\end{equation}
This right-hand side is 0 if
$c>-2\nu_{\Theta_\sigma}(x) \inf_{\lambda\ge 0}\log\bE[e^{-\theta_{\lambda,\sigma}(0)
/L(\lambda)}]$.
\end{proof}

By using the above path $r=(0=r_0,r_1,\dots,r_{N(r)}=nx)$, 
\begin{align*}
 &\frac{a_\lambda (0,nx,\omega ,\sigma )}{nL(\lambda)}\\
 &\quad\leq -\frac{1}{nL(\lambda)} \log E_\omega^x \Biggl[
 \exp \Biggl\{ -\sum_{m=0}^{H^X(y)-1} 
 \theta_{\lambda,\sigma} (X_m) \Biggr\}
       \,\1{\{ (X_m)_{m=0}^{N(r)}=r \}} \Biggr] \\
 &\quad\leq 
 \frac{1}{n}T(0,nx,\theta_{\lambda, \sigma}/L(\lambda))+\frac{1}{n} 
 +\frac{1}{L(\lambda)} \sum_{m=0}^{N(r)-1} 
 \frac{-\log \omega (r_m,r_{m+1}-r_m)}{n}.
\end{align*}
Since the last sum is bounded with high probability, that is,
\begin{equation*}
 \lim_{n\to\infty}\bP\left(\sum_{m=0}^{N(r)-1} 
 \frac{-\log \omega (r_m,r_{m+1}-r_m)}{n}
 \le 2c_x\E[-\log \max_{|e|=1}\omega(0,e)]\right)=1
\end{equation*}
by Lemma~\ref{percolation:lem_4} and the weak law of large numbers, 
we reach the desired conclusion.
\end{proof}

To complete the proof of Proposition~\ref{prop:percolation_2},
pick an arbitrary $\epsilon>0$ and take $\lambda>0$ so large that
$|\nu_{\theta_{\lambda,\sigma}/L(\lambda)}(x)-\nu_{\Theta_{\sigma}}(x)|<\epsilon$
and Lemma~\ref{lem:percolation_3} hold. 
Then we know that the events
\begin{align*}
&\left\{(\omega,\sigma):\left|\frac{\alpha_\lambda (x)}{L(\lambda)}-
\frac{a_\lambda (0,nx,\omega ,\sigma )}{nL(\lambda)}\right|
<\epsilon\right\},\\
&\left\{(\omega,\sigma):
\left|\frac{a_\lambda (0,nx,\omega ,\sigma )}{nL(\lambda)}
-\frac{1}{n}T(0,nx,\theta_{\lambda,\sigma}/L(\lambda))\right|
<\epsilon\right\},\\
&\left\{\sigma:\left|\frac{1}{n}T(0,nx,\theta_{\lambda,\sigma}/L(\lambda))
-\nu_{\theta_{\lambda,\sigma}/L(\lambda)}(x)\right|
<\epsilon\right\}
\end{align*}
have probability tending to 1 as $n\to\infty$.
In particular, we can find $(\omega,\sigma)$ belonging to all the
events above and substituting it into~\eqref{eq: percolation_4-1}, 
we obtain 
\begin{equation*}
\left|\frac{\alpha_\lambda (x)}{L(\lambda)}-\nu_{\Theta_\sigma} (x)\right|
<4\epsilon.
\end{equation*}
\end{proof}
\section{Some properties of the rate function}
\subsection{Asymptotics of the rate function}\label{rate function}
In this section, we discuss the asymptotics of the rate function as 
$x\to\infty$ and $x\to 0$ in some special cases. 

We start with the case $x\to\infty$.
Let 
\begin{equation*}
\lambda^*(x)=\inf\left\{\lambda\ge 0: 
L(\lambda)\nu_{\Theta_\sigma}(x)-\lambda = 
\sup_{\lambda\ge 0}(L(\lambda)\nu_{\Theta_\sigma}(x)-\lambda)\right\} 
\end{equation*}
with the convention $\inf\emptyset=\infty$.
\begin{proposition}\label{prop:rate_function_1}
Suppose that the same assumption as in 
Proposition~\ref{prop:percolation_2} holds. 
In addition, assume that $\lambda^*(x)<\infty$ for any $x\in\R^d$
and 
\begin{equation}
\lim_{\ell\to\infty}L(\lambda^*(\ell x))\nu_{\Theta_\sigma}(\ell x)/
\lambda^*(\ell x)>1. \label{eq:rate_function_1}
\end{equation}
Then for any $x\in\R^d$,
\begin{equation}
 I(\ell x)=\sup_{\lambda\ge 0}(L(\lambda)\ell \nu_{\Theta_\sigma}(x)-\lambda)
(1+o(1))\label{eq:rate_function_2}
\end{equation}
as $\ell\to\infty$.
\end{proposition}
\begin{proof}
Note that $\lambda^*(\ell x)\to\infty$ as $\ell \to\infty$. 
On the other hand, we know from Proposition~\ref{prop:percolation_2} 
that 
$\alpha_{\lambda}(x)= L(\lambda)\nu_{\Theta_\sigma}(x)(1+o(1))$ 
as $\lambda\to\infty$.
Combining these two facts and using \eqref{eq:rate_function_1}, 
we obtain
\begin{equation*}
\begin{split}
I(\ell x)&\ge
\alpha_{\lambda^*(\ell x)}(\ell x)-\lambda^*(\ell x)\\
&=L(\lambda^*(\ell x))\ell \nu_{\Theta_\sigma}(x)(1+o(1)))
-\lambda^*(\ell x)\\
&=(L(\lambda^*(\ell x))\ell \nu_{\Theta_\sigma}(x)-\lambda^*(\ell x))(1+o(1))
\end{split}
\end{equation*}
as $\ell\to\infty$. This proves the lower bound 
in~\eqref{eq:rate_function_2}. 
To prove the upper bound, fix any $\epsilon>0$. 
Then by the same reasoning as above, 
for sufficiently large $\ell$, we have 
\begin{equation*}
\sup_{\lambda\ge 0}
(\alpha_{\lambda}(\ell x)-\lambda)
\le \sup_{\lambda\ge 0}
((1+\epsilon)L(\lambda)\ell \nu_{\Theta_\sigma}(x)-\lambda)
\end{equation*}
and the right hand side is bounded from above by 
$(1+2\epsilon)\sup_{\lambda\ge 0}(\alpha_{\lambda}(\ell x)-\lambda)$.
\end{proof}

Using the above proposition, one can see that 
in the situation of Example~\ref{ex:percolation_0}-(2),
\begin{equation*}
I(\ell x)\sim \ell \nu_{\Theta_\sigma}(x)(\log(\ell \nu_{\Theta_\sigma}(x))-1) 
\textrm{ as }\ell \to \infty
\end{equation*}
and in that of Example~\ref{ex:percolation_0}-(4), 
\begin{equation*}
I(\ell x)\sim \frac{1}{1+\gamma}
\left(\frac{\gamma}{1+\gamma}\right)^{\gamma}
\left(\ell \nu_{\Theta_\sigma}(x)\right)^{1+\gamma}\textrm{ as }
\ell \to \infty.
\end{equation*}
Example~\ref{ex:percolation_0}-(1) does not fall within the 
scope of the above proposition but it is easy to see 
that $I(x)=\infty$ as soon as $\nu_{\Theta_\sigma}(x)>1$. 

Let us turn to the case $x\to 0$. We only consider the {\em simple 
random walk with random holding times}, i.e.,~$\omega(x,e)=\frac{1}{2d}$
for all $x\in\Z^d$ and $|e|=1$. This is of course very restrictive
but it seems rather unreasonable to expect a unified result under 
the general setting as large deviations of RWRE exhibit rich 
phenomena. For example, if a nestling RWRE satisfies the law of 
large numbers with nonzero speed $v$, then the rate function is zero
on the line segment connecting the origin and $v$. 

\begin{proposition}
Assume $\omega(x,e)=\frac{1}{2d}$ for all $x\in\Z^d$ and $|e|=1$
for $\P$ almost every $\omega$. 
Then
\begin{equation}
I(\ell x)=\frac{d}{2}\bE\left[\textstyle{\int_0^{\infty}s 
\sigma_0(ds)}\right]\ell^2|x|^2(1+o(1))\textrm{ as }\ell\to 0.
\label{eq:rate_function_3}
\end{equation}
\end{proposition}
\begin{proof}
Our $\alpha_{\lambda}$ is nothing but the quenched Lyapunov 
exponent of Green's function with the random potential 
$\theta_{\sigma,\lambda}$ (cf.~\cite{Zer98a}). 
It follows from our assumption 
$\int_0^{\infty}s \sigma_0(ds)\in L^d(\bP)$ that 
\begin{equation*}
\lambda^{-1}\theta_{\sigma,\lambda}(0)\to 
\int_0^{\infty}s \sigma_0(ds)\textrm{ as }\lambda\to 0
\end{equation*}
$\bP$-a.s.~and in $L^1(\bP)$. 
This verifies the assumption in Theorem~4 of~\cite{KMZ11}, 
which tells us that 
\begin{equation}
 \alpha_{\lambda}(x)
 =\sqrt{2d\lambda\bE\left[\textstyle{\int_0^{\infty}s \sigma_0(ds)}\right]}|x|
(1+o(1))\textrm{ as }\lambda\to 0.
\end{equation}
From this asymptotics, one can deduce~\eqref{eq:rate_function_3}
by the same way as for Proposition~\ref{prop:rate_function_1}.
\end{proof}

\subsection{Dependence of the rate function on the law of
the holding times} 
In this section, we discuss how the rate function is affected 
by the randomness of the holding times through simple examples.
More precisely, we establish some comparisons of the rate functions 
for different laws of holding times. 
We add superscripts to the rate functions and Lyapunov exponents, 
as in $I^\sigma$ and $\alpha_\lambda^\sigma$, to indicate their
dependence on the law of the holding times.

Let us begin with the comparison of the rate functions for 
a given law $\sigma=(\sigma_z)_{z\in \Z^d}$ of holding times with that 
for an {{averaged}} version 
$\bar{\sigma}=(\delta_{\int_0^\infty s \sigma_z(d s)})_{z\in\Z^d}$.
Under the latter law, the holding times are spatially inhomogeneous but deterministic
on each site. 
In this case, applying Jensen's inequality to the 
$E_\sigma^{\rm{HT}}$-expectation 
in the definition of $e_\lambda$, one immediately finds
$\alpha_\lambda^\sigma\le \alpha_\lambda^{\bar{\sigma}}$ 
and hence $I^\sigma\le I^{\bar{\sigma}}$.
Roughly speaking, this reflects the fact that random (rather than deterministic) holding 
times make it easier to realize rare events.

One can also consider another natural {averaged} version 
$\tilde\sigma_z(\cdot)=\int \sigma_z(\cdot) \bP(d\sigma)$. 
Then the holding times are random on each site but their joint law is 
spatially homogeneous. In this case, since
$e^{-\theta_{\lambda,\tilde\sigma}(z)}=\bE[e^{-\theta_{\lambda,\sigma}(z)}]$,
Jensen's inequality implies
\begin{equation*}
\bE\left[a_\lambda (x,y,\omega,\sigma )\right] 
\ge -\log \bE\left[e_\lambda (x,y,\omega ,\sigma )\right]
= a_\lambda (x,y,\omega,\tilde\sigma ).
\end{equation*}
From this and Theorem~\ref{thm:intro_1}, it follows that 
$\alpha_\lambda^\sigma\ge \alpha_\lambda^{\tilde{\sigma}}$ 
and hence $I^\sigma\ge I^{\tilde{\sigma}}$.
This means that the spatial inhomogeneity makes it difficult
for the random walk to reach any given remote point. 
Note also that $\alpha_\lambda^{\tilde{\sigma}}$ is related to
the annealed Lyapunov exponent studied in~\cite{Flu07}. 

Finally, we consider the case where $(\sigma_z^1)_{z\in\Z^d}$ and 
$(\sigma_z^2)_{z\in\Z^d}$ 
are collections of exponential distributions with i.i.d.~random rates 
$(r_z^1)_{z\in\Z^d}$ and $(r_z^2)_{z\in\Z^d}$. 
The exponential holding time is of special
importance since the RWREHT (when quenched) is a Markov process only in this case.
It is natural to expect that a certain variability of holding times 
plays a key role as above and we use the following notion:
\begin{definition}
Let $\mu$ and $\nu$ be probability distributions on $\R$. 
We say $\mu$ is more variable than $\nu$ if for every
concave increasing function $h:\R\to \R$, 
$\int h(x) \mu(dx)\le \int h(x) \nu(dx)$.
\end{definition}
By a direct computation and \eqref{eq:lyapunov_6}, we have
\begin{equation*}
 -a_\lambda(x,y,\omega,\sigma^i)
 =\log E_\omega^x \Biggl[  \Biggl(
\prod_{m=0}^{H^X(y)-1}\frac{r^i_{X_m}}{r^i_{X_m}+\lambda}\Biggr)
 \1{\{ H^X(y)<\infty \}} \Biggr].
\end{equation*}
This is a concave increasing function of $r_z^i$ for each $z\in\Z^d$.
Thus if the law of $r^1_0$ is more variable than that of $r^2_0$, 
it follows that
\begin{equation*}
 \bE[a_\lambda(x,y,\omega,\sigma^1)]\ge 
 \bE[a_\lambda(x,y,\omega,\sigma^2)]
\end{equation*}
by using Proposition~9.5.4 in~\cite{Ros96} just the same way as in 
Proposition~4 in~\cite{Zer98a}. 
As a result, $\alpha_\lambda^{\sigma_1} \ge \alpha_\lambda^{\sigma_2}$
and hence $I^{\sigma_1}\ge I^{\sigma_2}$, i.e., 
when the spatial inhomogeneity of the law of the holding times increases,
so does the rate function.

\subparagraph{acknowledgements.}
The authors would like to express their profound gratitude to the reviewer
for the very careful reading of the manuscript and also suggesting the 
consideration made in Subsection 6.2. 

\newcommand{\noop}[1]{}

\end{document}